\documentclass[11pt]{article}
\usepackage{mathrsfs}
\usepackage{amscd}
\usepackage{amsmath,amsfonts,amssymb,amscd}
\usepackage{indentfirst,graphics,epsfig,psfrag}
\input{epsf}
\usepackage{ifpdf}
\usepackage{enumerate}
\usepackage{appendix}
\usepackage{enumerate}
\usepackage{lineno}
\usepackage{color}
\makeatletter \@addtoreset{figure}{section} \makeatother
\makeatletter
\long\def\@makecaption#1#2{%
   \vskip 10\p@
   \setbox\@tempboxa\hbox{{#1}\ \ #2}%
   \ifdim \wd\@tempboxa >\hsize

       {#1}\ \ #2\par
   \else
       \hbox to\hsize{\hfil\box\@tempboxa\hfil}%
   \fi}
\makeatother

\newtheorem{thm}{Theorem}
\newtheorem{cor}{Corollary}
\newtheorem{lem}{Lemma}

\newtheorem{obs}{Observation}
\newtheorem{pro}{Proposition}

\newcommand{\qed}{{\hfill\rule{3pt}{7pt}}}

\def\qed{\hfill \rule{4pt}{7pt}}

\begin{document}
\title{\textbf{On the pedant tree-connectivity
of graphs}
\footnote{Supported by the National Science Foundation of China
(No. 11161037, 11101232, 11461054) and the Science Found of Qinghai Province
(No. 2014-ZJ-907).}}
\author{
\small Yaping Mao\footnote{E-mail: maoyaping@ymail.com}\\[0.3cm]
\small Department of Mathematics, Qinghai Normal\\
\small University, Xining, Qinghai 810008, China}
\date{}
\maketitle

\begin{abstract}
The concept of pedant tree-connectivity was introduced by Hager \cite{Hager} in
1985. For a graph $G=(V,E)$ and a set $S\subseteq V(G)$ of at least
two vertices, \emph{an $S$-Steiner tree} or \emph{a Steiner tree
connecting $S$} (or simply, \emph{an $S$-tree}) is a such subgraph
$T=(V',E')$ of $G$ that is a tree with $S\subseteq V'$. For an
$S$-Steiner tree, if the degree of each vertex in $S$ is equal to
one, then this tree is called a \emph{pedant $S$-Steiner tree}. Two
pedant $S$-Steiner trees $T$ and $T'$ are said to be
\emph{internally disjoint} if $E(T)\cap E(T')=\varnothing$ and
$V(T)\cap V(T')=S$. For $S\subseteq V(G)$ and $|S|\geq 2$, the
\emph{local pedant-tree connectivity} $\tau_G(S)$ is the maximum
number of internally disjoint pedant $S$-Steiner trees in $G$. For
an integer $k$ with $2\leq k\leq n$, \emph{$k$-pedant
tree-connectivity} is defined as
$\tau_k(G)=\min\{\tau_G(S)\,|\,S\subseteq V(G),|S|=k\}$. In this
paper, we first study the sharp bounds of pedant
tree-connectivity. Next, we obtain the exact value of a threshold graph, and give an
upper bound of the pedant-tree $k$-connectivity of a complete
multipartite graph. For a connected graph $G$, we show that $0\leq
\tau_k(G)\leq n-k$, and
graphs with $\tau_k(G)=n-k,n-k-1,n-k-2,0$ are characterized in this paper. In the end, we obtain the Nordhaus-Guddum type results for pedant tree-connectivity. \\[2mm]
\end{abstract}

\section{Introduction}

All graphs considered in this paper are undirected, finite and
simple. We refer to the book \cite{bondy} for graph theoretical
notation and terminology not described here. For a graph $G$, let
$V(G)$, $E(G)$ and $\delta(G)$ denote the set of vertices, the set
of edges and the minimum degree of $G$, respectively. Connectivity
is one of the most basic concepts of graph-theoretic subjects, both
in combinatorial sense and the algorithmic sense. As we know, the
classical connectivity has two equivalent definitions. The
\emph{connectivity} of $G$, written $\kappa(G)$, is the minimum
order of a vertex set $S\subseteq V(G)$ such that $G\setminus S$ is
disconnected or has only one vertex. We call this definition the
`cut' version definition of connectivity. A well-known theorem of
Whitney provides an equivalent definition of connectivity, which can
be called the `path' version definition of connectivity. For any two
distinct vertices $x$ and $y$ in $G$, the \emph{local connectivity}
$\kappa_{G}(x,y)$ is the maximum number of internally disjoint paths
connecting $x$ and $y$. Then
$\kappa(G)=\min\{\kappa_{G}(x,y)\,|\,x,y\in V(G),x\neq y\}$ is
defined to be the \emph{connectivity} of $G$. For connectivity,
Oellermann gave a survey paper on this subject; see
\cite{Oellermann2}.

Although there are many elegant and powerful results on connectivity
in Graph Theory, the basic notation of classical connectivity may
not be general enough to capture some computational settings. So
people want to generalize this concept. For the `cut' version
definition of connectivity, we find the above minimum vertex set
without regard the number of components of $G\setminus S$. Two
graphs with the same connectivity may have differing degrees of
vulnerability in the sense that the deletion of a vertex cut-set of
minimum cardinality from one graph may produce a graph with
considerably more components than in the case of the other graph.
For example, the star $K_{1,n}$ and the path $P_{n+1}\ (n\geq 3)$
are both trees of order $n+1$ and therefore connectivity $1$, but
the deletion of a cut-vertex from $K_{1,n}$ produces a graph with
$n$ components while the deletion of a cut-vertex from $P_{n+1}$
produces only two components. Chartrand et al. \cite{Chartrand1}
generalized the `cut' version definition of connectivity. For an
integer $k \ (k\geq 2)$ and a graph $G$ of order $n \ (n\geq k)$,
the \emph{$k$-connectivity} $\kappa'_k(G)$ is the smallest number of
vertices whose removal from $G$ of order $n \ (n\geq k)$ produces a
graph with at least $k$ components or a graph with fewer than $k$
vertices. Thus, for $k=2$, $\kappa'_2(G)=\kappa(G)$. For more
details about $k$-connectivity, we refer to \cite{Chartrand1, Day,
Oellermann2, Oellermann3}.

The generalized connectivity of a graph $G$, introduced by Hager \cite{Hager}, is
a natural generalization of the `path' version definition of
connectivity. For a graph $G=(V,E)$ and a set $S\subseteq V(G)$ of
at least two vertices, \emph{an $S$-Steiner tree} or \emph{a Steiner
tree connecting $S$} (or simply, \emph{an $S$-tree}) is a such
subgraph $T=(V',E')$ of $G$ that is a tree with $S\subseteq V'$.
Note that when $|S|=2$ a Steiner tree connecting $S$ is just a path
connecting the two vertices of $S$. Two Steiner trees $T$ and $T'$
connecting $S$ are said to be \emph{internally disjoint} if
$E(T)\cap E(T')=\varnothing$ and $V(T)\cap V(T')=S$. For $S\subseteq
V(G)$ and $|S|\geq 2$, the \emph{generalized local connectivity}
$\kappa_G(S)$ is the maximum number of internally disjoint trees
connecting $S$ in $G$, that is, we search for the maximum
cardinality of edge-disjoint trees which include $S$ and are vertex
disjoint with the exception of $S$. For an integer $k$ with $2\leq
k\leq n$, \emph{generalized $k$-connectivity} (or
\emph{$k$-tree-connectivity}) is defined as
$\kappa_k(G)=\min\{\kappa_G(S)\,|\,S\subseteq V(G),|S|=k\}$, that
is, $\kappa_k(G)$ is the minimum value of $\kappa_G(S)$ when $S$
runs over all $k$-subsets of $V(G)$. Clearly, when $|S|=2$,
$\kappa_2(G)$ is nothing new but the connectivity $\kappa(G)$ of
$G$, that is, $\kappa_2(G)=\kappa(G)$, which is the reason why one
addresses $\kappa_k(G)$ as the generalized connectivity of $G$. By
convention, for a connected graph $G$ with less than $k$ vertices,
we set $\kappa_k(G)=1$. Set $\kappa_k(G)=0$ when $G$ is
disconnected. Note that the generalized $k$-connectivity and
$k$-connectivity of a graph are indeed different. Take for example,
the graph $H_1$ obtained from a triangle with vertex set
$\{v_1,v_2,v_3\}$ by adding three new vertices $u_1,u_2,u_3$ and
joining $v_i$ to $u_i$ by an edge for $1 \leq i\leq 3$. Then
$\kappa_3(H_1)=1$ but $\kappa'_3(H_1)=2$. There are many results on
the generalized connectivity, see \cite{Chartrand2, LLMS,
LLSun, LL, LLZ, LM1, LM2, LM3, LM4, LMS, Okamoto}.

The concept of pedant-tree connectivity \cite{Hager} was introduced
by Hager in 1985, which is specialization of generalized
connectivity (or \emph{$k$-tree-connectivity}) but a generalization of classical connectivity. For an $S$-Steiner
tree, if the degree of each vertex in $S$ is equal to one, then this
tree is called a \emph{pedant $S$-Steiner tree}. Two pedant
$S$-Steiner trees $T$ and $T'$ are said to be \emph{internally
disjoint} if $E(T)\cap E(T')=\varnothing$ and $V(T)\cap V(T')=S$.
For $S\subseteq V(G)$ and $|S|\geq 2$, the \emph{local-pedant
tree connectivity} $\tau_G(S)$ is the maximum number of internally
disjoint pedant $S$-Steiner trees in $G$. For an integer $k$ with
$2\leq k\leq n$, \emph{pedant-tree $k$-connectivity} is defined as
$\tau_k(G)=\min\{\tau_G(S)\,|\,S\subseteq V(G),|S|=k\}$. Set
$\kappa_k(G)=0$ when $G$ is disconnected.

In {\upshape\cite{Hager}}, Hager derived the following results.

\begin{lem}{\upshape\cite{Hager}}\label{lem1}
Let $G$ be a graph. If $\tau_k(G)\geq \ell$, then $\delta(G)\geq
k+\ell-1$.
\end{lem}

\begin{lem}{\upshape\cite{Hager}}\label{lem2}
Let $G$ be a graph. If $\tau_k(G)\geq \ell$, then $\kappa(G)\geq
k+\ell-2$.
\end{lem}

Li et al. {\upshape\cite{LLZ}} obtained the following result.

\begin{lem}{\upshape\cite{LLZ}}\label{lem3}
Let $G$ be a connected graph with minimum degree $\delta$. If there are two adjacent vertices of degree $\delta$, then
$\kappa(G)\leq \delta(G)-1$.
\end{lem}

Obviously, the generalized $k$-connectivity (or
\emph{$k$-tree-connectivity}) and pedant-tree $k$-connectivity of a
graph are indeed different. For example, let $H=W_n$ be a wheel of
order $n$. From Lemma \ref{lem1}, we have $\tau_3(H)\leq 1$. One
can check that for any $S\subseteq V(H)$ with $|S|=3$,
$\tau_{H}(S)\geq 1$. Therefore, $\tau_3(H)=1$. From Lemma
\ref{lem3}, we have $\kappa_3(H)\leq \delta(H)-1=3-1=2$. One can
check that for any $S\subseteq V(G)$ with $|S|=3$, $\kappa_{H}(S)\geq
2$. Therefore, $\kappa_3(H)=2$.

In addition to being a natural combinatorial measure, both the
pendant tree-connectivity and the generalized connectivity can be
motivated by its interesting interpretation in practice. For
example, suppose that $G$ represents a network. If one considers to
connect a pair of vertices of $G$, then a path is used to connect
them. However, if one wants to connect a set $S$ of vertices of $G$
with $|S|\geq 3$, then a tree has to be used to connect them. This
kind of tree with minimum order for connecting a set of vertices is
usually called a Steiner tree, and popularly used in the physical
design of VLSI (see \cite{Grotschel1, Grotschel2, Sherwani}) and
computer communication networks (see \cite{Du}). Usually, one wants
to consider how tough a network can be, for the connection of a set
of vertices. Then, the number of totally independent ways to connect
them is a measure for this purpose. The generalized $k$-connectivity
can serve for measuring the capability of a network $G$ to connect
any $k$ vertices in $G$.

This paper is organized as follows. In Section $2$, we study sharp bounds of pedant tree-connectivity. In Section $3$, we obtain the exact value of a threshold graph, and give an upper bound of the pedant-tree $k$-connectivity of a complete
multipartite graph. For a connected graph $G$, we show that $0\leq
\tau_k(G)\leq n-k$, and graphs with $\tau_k(G)=n-k,n-k-1,n-k-2$ are characterized in Section $4$.
By Fan Lemma, graphs with $\tau_k(G)=0$ are characterized in Section $5$.

Let $\mathcal {G}(n)$ denote the class of simple graphs of order $n$
and $\mathcal {G}(n,m)$ the subclass of $\mathcal {G}(n)$ having $m$
edges. Give a graph theoretic parameter $f(G)$ and a positive
integer $n$, the \emph{Nordhaus-Gaddum(\textbf{N-G}) Problem} is to
determine sharp bounds for: $(1)$ $f(G)+f(\bar{G})$ and $(2)$
$f(G)\cdot f(\bar{G})$, as $G$ ranges over the class $\mathcal
{G}(n)$, and characterize the extremal graphs. The Nordhaus-Gaddum
type relations have received wide investigations. Recently,
Aouchiche and Hansen published a survey paper on this subject, see
\cite{Aouchiche}. In Section $6$, we study the Nordhaus-Gaddum problem for pedant tree-connectivity.

\section{Sharp bounds of pedant tree-connectivity}

In {\upshape\cite{Hager}}, Hager obtained the exact value of the pedant-tree $k$-connectivity of a complete graph.

\begin{lem}{\upshape\cite{Hager}}\label{lem4}
Let $k,n$ be two integers with $3\leq k\leq n$, and let $K_n$ be a
complete graph of order $n$. Then
$$
\tau_k(K_n)=n-k.
$$
\end{lem}

Let $S$ be a set of $k$ vertices of a connected graph $G$, and let
$\mathcal{T}$ be a set of internally disjoint pedant $S$-Steiner trees.
Then the following observation is immediate.

\begin{obs}\label{obs1}
Let $k,n$ be two integers with $3\leq k\leq n$. Let $G$ be a
connected graph of order $n$, and let $S\subseteq V(G)$ with $|S|=k$.
For each $T\in \mathcal{T}$,
$$
|E(T)\cap E_G[S,\bar{S}]|\geq k,
$$
where $\bar{S}=V(G)\setminus S$.
\end{obs}

From the above result, we can derive an upper bound of $k$-pedant
tree-connectivity.

\begin{thm}\label{th1}
For any graph $G$ with order at least $k$,
$$
\tau_k(G)\leq \min_{S\subseteq V(G),
|S|=k}\Big\lfloor\frac{1}{k}|E_G[S,\bar{S}]|\Big\rfloor,
$$
where
$S\subseteq V(G)$ with $|S|=k$, and $\bar{S}=V(G)\setminus S$.
Moreover, the bound is sharp.
\end{thm}
\begin{pf}
For any $S\subseteq V(G)$ with $|S|=k$, it suffices to
$\tau(S)=|\mathcal{T}|\leq \frac{1}{k}|E_G[S,\bar{S}]|$. From
Observation \ref{obs1}, for each tree $T\in \mathcal{T}$, we have
$$
|E(T)\cap E_G[S,\bar{S}]|\geq k.
$$
Therefore,
$$
k|\mathcal{T}|\leq |E_G[S,\bar{S}]|.
$$
Since $|\mathcal{T}|$ is an integer, we have
$$
\tau(S)=|\mathcal{T}|\leq
\Big\lfloor\frac{1}{k}|E_G[S,\bar{S}]|\Big\rfloor.
$$
From the arbitrariness of $S$, we have
$$
\tau_k(G)=\min_{S\subseteq V(G), |S|=k} \tau(S)=\min_{S\subseteq
V(G), |S|=k}|\mathcal{T}|\leq \min_{S\subseteq V(G),
|S|=k}\Big\lfloor\frac{1}{k}|E_G[S,\bar{S}]|\Big\rfloor.
$$

To show the bound is sharp, we consider the graph $G=K_n$. From
Lemma \ref{lem4}, we have $\tau_k(G)=n-k$. Thus
$|E_{K_n}[S,\bar{S}]|=k(n-k)$ and hence
$\frac{1}{k}|E_{K_n}[S,\bar{S}]|=n-k$. Therefore,
$$
\min_{S\subseteq V(G),
|S|=k}\Big\lfloor\frac{1}{k}|E_G[S,\bar{S}]|\Big\rfloor=n-k=\tau_k(G).
$$
So the bound of this theorem is sharp. For $k=3$, the graph
$H=W_n$ is a sharp example. As we know, $\tau_3(W_n)=1$. Observe that the
graph $H=W_n$ is obtained from a cycle $C=v_1v_2\cdots v_{n-1}v_1$ by
adding a vertex $v_n$ and the edges $v_iv_n \ (1\leq i\leq n-1)$.
Choose $S=\{v_2,v_3,v_4\}$. Then
$E_{H}[S,\bar{S}]=\{v_1v_2,v_2v_n,v_3v_n,v_4v_n,v_4v_5\}$ and
hence we have
$\lfloor\frac{1}{3}|E_{H}[S,\bar{S}]|\rfloor=\lfloor\frac{1}{3}\cdot
5\rfloor=1=\tau_3(H)$. \qed
\end{pf}\\

For any connected graph $G$ of order $n$, from Lemma \ref{lem4}, we
have $\tau_k(G)\leq \tau_k(K_n)=n-k$. The following upper and lower bounds for
$\tau_k(G)$ can be easily seen.

\begin{pro}\label{pro1}
Let $k,n$ be two integers with $3\leq k\leq n$, and let $G$ be a
graph. Then
$$
0\leq \tau_k(G)\leq n-k.
$$
\end{pro}

For $k=n$, the following corollary is immediate.

\begin{cor}\label{cor1}
Let $G$ be a graph of order $n$. Then $\tau_n(G)=0$ if and only if
$G$ is a graph of order $n$.
\end{cor}

\section{Pedant-tree connectivity of some special graphs}

In \cite{Hager}, Hager obtained the exact value of pedant tree-connectivity of a complete bipartite graph.

\begin{lem}{\upshape\cite{Hager}}\label{lem5}
Let $K_{r,s}$ be a complete bipartite graph with $r+s$ vertices.
Then
$$
\tau_k(K_{r,s})=\max\{\min\{r-k+1, s-k+1\},0\}.
$$
\end{lem}

For complete multipartite graphs, we obtain the following result.

\begin{thm}\label{th2}
Let $K_{n_1,n_2,\cdots,n_t}$ be a complete $t$-partite graph with
$n_1\leq n_2\leq \cdots \leq n_t$.

$(1)$ If $k\geq t$, then
$$
\tau_k(K_{n_1,n_2,\cdots,n_t})\leq
\left\lfloor\frac{\sum_{i=1}^tn_i-k}{2}\right\rfloor.
$$

$(2)$ If $k<t$, then
$$
\tau_k(K_{n_1,n_2,\cdots,n_t})\leq \sum_{i=k+1}^tn_i+
\left\lfloor\frac{\sum_{i=1}^kn_i-k}{2}\right\rfloor.
$$
Moreover, the upper bounds are sharp.
\end{thm}
\begin{pf}
Set $G=K_{n_1,n_2,\cdots,n_t}$. Let $V_1,V_2,\cdots,V_t$ be the parts of complete $t$-partite graph $G$, and set $|V_i|=n_i \ (1\leq i\leq t)$.

$(1)$ Suppose $k\geq t$. Choose $S\subseteq V(G)$ and $|S|=k$ such that $S\cap V_i\neq \varnothing  \ (1\leq i\leq t)$. Observe that any pedant $S$-Steiner tree must muse at least two vertices in $\bar{S}=V(G)\setminus S$. Therefore,
$$
\tau_k(G)\leq \tau(S)\leq
\left\lfloor\frac{\sum_{i=1}^tn_i-k}{2}\right\rfloor.
$$

$(2)$ Suppose $k<t$. Choose $S\subseteq V(G)$ and $|S|=k$ such that $S\cap V_i\neq \varnothing  \ (1\leq i\leq k)$. Then $S\cap V_i=\varnothing  \ (k+1\leq i\leq t)$. Observe that any pedant $S$-Steiner tree must muse at least one vertex in $\bar{S}=V(G)\setminus S$. Since $n_1\leq n_2\leq \cdots \leq n_t$, it follows that there are at most $\sum_{i=k+1}^tn_i$ pedant $S$-Steiner trees in $G$. Any other pedant $S$-Steiner tree must use at least two vertices in $\bar{S}$. Therefore,
$$
\tau_k(G)\leq \tau(S)\leq
\sum_{i=k+1}^tn_i+
\left\lfloor\frac{\sum_{i=1}^kn_i-k}{2}\right\rfloor.
$$
The proof is now complete.\qed
\end{pf}

To show the sharpness of the bound in Theorem \ref{th2}, we consider the following example.\\

\noindent{\bf Example 1:}  Let $G=K_{r,r,r}$ be complete $3$-part graph where $r=n_1=n_2=n_3$. Suppose $k=3$. From the above theorem, $\tau_3(G)\leq \lfloor\frac{3r-3}{2}\rfloor$. To show $\tau_3(G)\geq  \lfloor\frac{3r-3}{2}\rfloor$, it suffices to prove that $\tau(S)\geq  \lfloor\frac{3r-3}{2}\rfloor$ for any $S\subseteq V(G)$ and $|S|=3$. Let $U=\{u_1,u_2,\cdots,u_r\}$, $V=\{v_1,v_2,\cdots,v_r\}$ and $W=\{w_1,w_2,\cdots,w_r\}$ be the three parts of $G$. Suppose $|S\cap U|=3$ or $|S\cap V|=3$ or $|S\cap W|=3$. Without loss of generality, let $|S\cap U|=3$. Then the trees $T_i$ induced by the
edges in $\{v_iu_1,v_iu_2,v_iu_3\} \ (1\leq i\leq r)$ and the trees $T_j'$ induced by the
edges in $\{w_ju_1,w_ju_2,w_ju_3\} \ (1\leq j\leq r)$ are
$2r>\lfloor\frac{3r-3}{2}\rfloor$ internally disjoint pedant $S$-Steiner trees, which implies $\tau(S)\geq \lfloor\frac{3r-3}{2}\rfloor$, as desired.
Suppose $|S\cap U|=2$ or $|S\cap V|=2$ or $|S\cap W|=2$. Without loss of generality, let $|S\cap U|=2$. Then $|S\cap V|=1$ or $|S\cap W|=1$. Without loss of generality, let $S=\{u_1,u_2,v_1\}$. Then the trees $T_i$ induced by the
edges in $\{w_iu_1,w_iu_2,w_iu_3\} \ (1\leq i\leq r)$ and the trees $T_j'$ induced by the
edges in $\{v_ju_1,v_ju_2,u_jv_1,u_jv_j\} \ (3\leq j\leq r)$ are
$r+\lfloor\frac{r-2}{2}\rfloor\geq \lfloor\frac{3r-3}{2}\rfloor$ internally disjoint pedant $S$-Steiner trees, and hence $\tau(S)\geq \lfloor\frac{3r-3}{2}\rfloor$, as desired.
Suppose $|S\cap U|=|S\cap V|=|S\cap W|=1$. Without loss of generality, let $S=\{u_r,v_r,w_r\}$. Then the trees $T_i$ induced by the
edges in $\{v_iu_r,u_iv_r,u_iv_i,w_iw_r,v_iw_r\}$ and the trees $T_i'$ induced by the
edges in $\{v_{i+1}u_r,v_{i+1}w_r,v_{i+1}w_i,w_iv_r\}$
and the trees $T_i''$ induced by the
edges in $\{u_{i+1}w_r,u_{i+1}w_{i+1},w_{i+1}v_r,w_{i+1}u_r\}$
are
$\lfloor\frac{3r-3}{2}\rfloor$ internally disjoint pedant $S$-Steiner trees, and hence $\tau(S)\geq \lfloor\frac{3r-3}{2}\rfloor$, where $i=2\ell-1$ and $1\leq \ell\leq \lfloor\frac{r-1}{2}\rfloor$. From the above argument, $\tau_3(G)=\lfloor\frac{3r-3}{2}\rfloor$, which implies that the bound of $(1)$ of Theorem \ref{th2} is sharp. For the bound of $(2)$ of Theorem \ref{th2}, one can check that the complete $4$-partite graph $K_{r,r,r,r}$ is a sharp example.\\

A graph $G$ is a \emph{threshold graph}, if there exists a weight
function $w: V(G)\rightarrow R$ and a real constant $t$ such that
two vertices $u,v\in V(G)$ are adjacent if and only if
$w(u)+w(v)\geq t$.

The following observation is easy to make from the definition of a
threshold graph.

\begin{obs}\label{obs2}
Let $G([n],E)$ be a threshold graph with a weight function $w:
V(G)\rightarrow R$. Let the vertices be labelled so that $w(1)\geq w(2)\geq
\cdots \geq w(n)$. Then

$(a)$ $d_1\geq d_2\geq \cdots \geq d_n$, where $d_i$ is the degree of vertex
$i$.

$(b)$ $I=\{i\in V(G): d_i\leq i-1\}$ is a maximum independent set of $G$
and $G\setminus I$ is a clique in $G$.

$(c)$ $N(i)=\{1,2,\cdots,d_i\}$ for every $i\in I$. Thus, the
neighborhoods of vertices in $I$ form a linear order under set
inclusion. Furthermore, if $G$ is connected, then every vertex in $G$ is
adjacent to $1$.
\end{obs}

Now, we are in a position to give our result.

\begin{thm}\label{th3}
Let $G$ be a threshold graph with $\delta(G)=\ell$. Then
$$
\tau_k(G)=\left\{
\begin{array}{ll}
0,&if~k>\ell;\\
\ell-k+1,&if~k\leq \ell.
\end{array}
\right.
$$
\end{thm}
\begin{pf}
Let $C_r$ and $I_{n-r}$ denote the clique and the maximum
independent set of $G$, respectively. Set
$V(C_r)=\{u_1,u_2,\cdots,u_r\}$ and
$V(I_{n-r})=\{v_1,v_2,\cdots,v_{n-r}\}$. Since $\delta(G)=\ell$, it
follows that $v_ju_i\in E(G)$ for each $u_i \ (1\leq i\leq \ell)$
and each $v_j \ (1\leq j\leq n-r)$.

Suppose $k\geq \ell+1$. Choose
$S=\{u_1,u_2,\cdots,u_{k-1},v_{n-r}\}$. Observe that $v_{n-r}u_i\in
E(G) \ (1\leq i\leq \ell)$ and $v_{n-r}u_i\notin E(G) \ (\ell+1 \leq
i\leq r)$. Then any pedant $S$-Steiner tree must occupy
some edge $v_{n-r}u_j\in E(G) \ (1\leq j\leq \ell)$. Since $u_j\in
S \ (1\leq j\leq \ell)$, it follows that the degree of $u_j$ in this tree is at least two, a
contradiction. So $\tau_k(G)\leq 0$. Combining this with Proposition
\ref{pro1}, we have $\tau_k(G)=0$.

Suppose $k\leq \ell$. Choose $S=\{u_1,u_2,\cdots,u_{k-1},v_{n-r}\}$.
Note that $v_{n-r}u_i\in E(G) \ (1\leq i\leq \ell)$ and
$v_{n-r}u_i\notin E(G) \ (\ell+1 \leq i\leq r)$. Then any
pedant $S$-Steiner tree $T$ must occupy some edge
$v_{n-r}u_j\in E(G) \ (1\leq j\leq \ell)$. If $v_{n-r}u_j\in E(T) \
(1\leq j\leq k-1)$, then the degree of $u_j$ in $T$ is at least two, a
contradiction. We now assume $v_{n-r}u_j\in E(G) \ (k\leq j\leq
\ell)$. Because a pedant $S$-Steiner tree must occupy at
least one edge in $\{v_{n-r}u_j\,|\,k\leq j\leq \ell\}$, we have at most
$\ell-k+1$ pedant $S$-Steiner trees in $G$, which implies
$\tau_k(G)\leq \ell-k+1$.

From the definition of $\tau_k(G)$,
it suffices to show that $\tau(S)\geq \ell-k+1$ for any $S\subseteq
V(G)$ and $|S|=k$. Suppose $S\subseteq V(C_r)$. Set
$S=\{u_{i_1},u_{i_2},\cdots,u_{i_k}\}$. From Lemma \ref{lem4}, $\tau_k(C_r)= r-k\geq \ell -k$, and hence there are $\ell-k$ pedant $S$-Steiner trees in $G$. These trees together with the tree induced by the edges in $\{u_{i_1}v_1,u_{i_2}v_1,\cdots,u_{i_k}v_1\}$ are $\ell-k+1$ pedant $S$-Steiner trees in $G$, as desired. Suppose $S\subseteq V(I_{n-r})$. Since $u_jv_{i_t}\in E(G) \ (1\leq j\leq \ell, 1\leq t\leq k)$, it follows that the trees induced by the edges in $\{v_{i_1}u_j,v_{i_2}u_j,\cdots,v_{i_k}u_j\}$ are $\ell$ pedant $S$-Steiner trees in $G$, as desired. Suppose $S\cap V(C_r)\neq \varnothing$ and $S\cap V(I_{n-r})\neq \varnothing$. Set $|S\cap V(C_r)|=p$ and $|S\cap V(I_{n-r})|=k-p$. Clearly, $1\leq p\leq k-1$ and $|S\cap \{u_{1},u_{2},\cdots,u_{\ell}\}|\leq p$. Then there exist at least $\ell-p$ vertices belonging to $V(G)\setminus S$. Choose $\ell-k+1$ vertices from them, say $u_{j_1},u_{j_2},\cdots,u_{j_{\ell-k+1}}$. Set
$S\cap V(C_r)=\{u_{i_1},u_{i_2},\cdots,u_{i_p}\}$. Without loss of generality, let $S\cap V(I_{n-r})=\{v_{1},v_{2},\cdots,u_{k-p}\}$. Then the trees induced by the edges in $\{u_{i_1}u_{j_i},u_{i_2}u_{j_i},\cdots,u_{i_p}u_{j_i}\}\cup \{v_{1}u_{j_i},v_{2}u_{j_i},\cdots,v_{k-p}u_{j_i}\} \ (1\leq i\leq \ell-k+1)$ are $\ell-k+1$ pedant $S$-Steiner trees in $G$, as desired. From the above argument, we know that $\tau_k(G)=\ell-k+1$. \qed
\end{pf}

\section{Graphs with large pedant tree-connectivity}

The graphs attaining the upper bound of Proposition \ref{pro1} can
be characterized now.

\begin{thm}\label{th4}
Let $k,n$ be two integers with $3\leq k\leq n$ and $n\geq 4$, and
let $G$ be a connected graph. Then $\tau_k(G)=n-k$ if and only if $G$ is
a complete graph of order $n$.
\end{thm}
\begin{pf}
Suppose $\tau_k(G)=n-k$. From Lemma \ref{lem1}, we have
$\delta(G)\geq k+(n-k)-1=n-1$. So $G$ is a complete graph of order
$n$. Conversely, suppose $G$ is a complete graph of order $n$. For
any $S\subseteq V(G)$ with $|S|=k$, we have $|V(G)\setminus S|=n-k$.
Let $S=\{u_1,u_2,\cdots,u_k\}$ and $V(G)\setminus
S=\{w_1,w_2,\cdots,w_{n-k}\}$. Then the trees $T_i$ induced by the
edges in $\{w_iu_1,w_iu_2,\cdots,w_iu_k\} \ (1\leq k\leq n-k)$ are
$n-k$ internally disjoint pedant $S$-Steiner tree, which
implies that $\tau(S)\geq n-k$. From the arbitrariness of $S$, we
have $\tau_k(G)\geq n-k$. Combining this with Proposition
\ref{pro1}, $\tau_k(G)=n-k$. \qed
\end{pf}\\

The graphs with $\tau_k(G)=n-k-1$ can also be characterized in the
following.

\begin{thm}\label{th5}
Let $k,n$ be two integers with $3\leq k\leq n$ and $n\geq 7$, and let $G$ be a
connected graph. Then $\tau_k(G)=n-k-1$ if and only if $\bar{G}=r K_2\cup (n-2r)K_1\ (r=1,2)$.
\end{thm}
\begin{pf}
Suppose $\tau_k(G)=n-k-1$. From Lemma \ref{lem1}, we have
$\delta(G)\geq k+(n-k-1)-1=n-2$. Therefore, $G$ is a graph obtained from a complete graph
of order $n$ by deleting a matching $M$ such that $0\leq |M|\leq
\lfloor\frac{n}{2}\rfloor$. Then we have the following claim.

{\bf Claim 1.} $1\leq |M|\leq
2$.

\noindent{\bf Proof of Claim 1.} Set $|M|=r$. Then $1\leq r\leq
\lfloor\frac{n}{2}\rfloor$. Assume $r\geq 3$. Without loss of
generality, let $M=\{u_iw_i\,|\,1\leq i\leq r\}$. Choose
$S=\{u_1,u_2,u_3,v_1,v_2,\cdots,v_{k-3}\}\subseteq V(G)$ where
$v_1,v_2,\cdots,v_{k-3}\in V(G)\setminus
\{u_1,u_2,u_3,w_1,w_2,w_3\}$. Thus $u_iw_i\notin E_G[S,\bar{S}]
(1\leq i\leq 3)$, where $\bar{S}=V(G)\setminus S$. From the definition of $\tau_k(G)$, any pedant $S$-Steiner tree must use at least one vertex of $\bar{S}$. Observe that any
pedant $S$-Steiner tree containing vertex $w_i \
(1\leq i\leq 3)$ must occupy at least two vertices of $\{w_1,w_2,w_3\}$ since $u_iw_i\notin E_G[S,\bar{S}] \ (1\leq i\leq 3)$. So
$w_1,w_2,w_3$ must belong to two pedant $S$-Steiner tree,
say $T_1,T_2$, which implies that these trees occupy at least
four vertices of $\bar{S}=V(G)\setminus S$. So there are at most $n-k-2$
internally disjoint pedant $S$-Steiner trees, a contradiction.

From Claim $1$, we have $1\leq |M|\leq
2$, as desired.

Conversely, suppose that $G$ is a connected graph of order $n$ satisfying
the conditions of this theorem. It suffices to show
$\tau_k(G)\geq n-k-1$, where $G=K_n\setminus M$ and $M$ is a
matching of size $1$ or $2$. In fact, we only need to prove
$\tau_k(G)\geq n-k-1$, where $G=K_n\setminus M$ and $M$ is a
matching of size $2$. From the definition
of $\tau_k(G)$, it suffices to show that $\tau(S)\geq n-k-1$ for any
$S\subseteq V(G)$ with $|S|=k$. Since $|M|=2$, it follows that
$0\leq |M\cap E_{K_n}[S,\bar{S}]|\leq 2$.

If $|M\cap E_{K_n}[S,\bar{S}]|=2$, then we set
$M=\{u_1v_1,u_2v_2\}$ and $V(G)\setminus
\{u_1,u_2,v_1,v_2\}=\{w_1,w_2,\cdots,w_{n-4}\}$. Without loss of generality, let $S=\{u_1,u_2,w_1,w_2,\cdots,w_{k-2}\}$.
Then $v_1,v_2\in \bar{S}$. Then the trees $T_j$ induced by the edges
in $\{w_ju_1,w_ju_2,w_jw_1, v_jw_2,$ $\cdots,w_ju_{k-2}\} \ (k-1\leq
j\leq n-4)$ together with the tree $T_1$ induced by the edges in
$\{v_1u_2,v_1w_1,v_1w_2,\cdots,v_1w_{k-2},$ $v_1v_2,v_2u_1\}$ form
$(n-4)-(k-2)+1=n-k-1$ pedant $S$-Steiner trees, which
implies $\tau(S)\geq n-k-1$.

If $|M\cap E_{K_n}[S,\bar{S}]|=1$, then we set $u_1v_1\in M\cap E_{K_n}[S,\bar{S}]$ and $V(G)\setminus
\{u_1,v_1\}=\{w_1,w_2,\cdots,w_{n-2}\}$.
Without loss of generality, let $S=\{u_1,w_1,w_2,\cdots,w_{k-1}\}$.
Then the trees $T_j$ induced
by the edges in $\{w_ju_1,w_jw_1,w_jw_2,\cdots,w_ju_k\} \ (k\leq
j\leq n-2)$ form $(n-4)-(k-1)+2=n-k-1$ pedant $S$-Steiner trees, and hence $\tau(S)\geq n-k-1$.

If $|M\cap E_{K_{n}}[S,\bar{S}]|=0$, then we let $V(G)=\{w_1,w_2,\cdots,w_{n}\}$. Without loss of generality, let $S=\{w_1,w_2,\cdots,w_k\}$. Then the
trees $T_j$ induced by the edges in $\{w_jw_1,w_jw_2,\cdots,
w_jw_k\} \ (k+1\leq j\leq n-4)$ together with the trees $T_i$
induced by the edges in $\{w_1u_i,w_2u_i,\cdots,w_ku_i\} \ (1\leq
i\leq 2)$ and the trees $T'_j$ induced by the edges in $\{w_1v_j,
w_2v_j,\cdots,w_kv_j\} \ (1\leq j\leq 2)$ form $(n-4)-k+2+2=n-k$
pedant $S$-Steiner trees, which implies that $\tau(S)\geq
n-k$.

From the above argument, we conclude that for any $S\subseteq V(G)$
with $|S|=k$, $\tau(S)\geq n-k-1$. From the arbitrariness of $S$, we
have $\tau_k(G)\geq n-k-1$. Combining this with Theorem \ref{th4},
$\tau_k(G)=n-k-1$. \qed
\end{pf}\\

For $k=3$, graphs with $\tau_k(G)=n-k-2$ are characterized in the following lemma, which is preparation of Theorem \ref{th6}.

\begin{lem}\label{lem6}
Let $G$ be a connected graph of order $n$.

$(1)$ For $k=3$, $\tau_k(G)=n-k-2$ if and only if $\bar{G}$ is a subgraph of one of the following graphs.
\begin{itemize}
\item $C_i\cup C_j\cup (n-i-j)K_1 \ \ (i=3,4, \ j=3,4)$;

\item $C_i\cup \lfloor \frac{n-i}{2} \rfloor K_2  \ \ (i=3,4)$;

\item $P_5\cup \lfloor \frac{n-5}{2} \rfloor K_2$;

\item $C_i\cup (n-i) K_1  \ \ (i=5,6,7)$.
\end{itemize}

$(2)$ For $k=4$, $\tau_k(G)=n-k-2$ if and only if $\bar{G}$ is a subgraph of one of the following graphs.
\begin{itemize}
\item $C_i\cup C_j\cup (n-i-j)K_1 \ \ (i=3,4, \ j=3,4)$;

\item $C_i\cup \lfloor \frac{n-i}{2} \rfloor K_2  \ \ (i=3,4)$;

\item $C_5\cup K_2\cup (n-7) K_1$;

\item $C_i\cup (n-i) K_1  \ \ (i=6,7)$.
\end{itemize}
\end{lem}

\begin{pf}
Suppose $\tau_k(G)=n-k-2 \ (k=3,4)$. We first give the proof of the following claim.

{\bf Claim 1.} For any $S\subseteq V(G)$ and $|S|=k$, if $|N_{\bar{G}}(S)|\geq 5$, then $\tau_k(G)\leq n-k-3$.

\noindent{\bf Proof of Claim 1:} Assume, to the contrary, that $\tau_k(G)\geq n-k-2$.
Choose five vertices in $N_{\bar{G}}(S)$, say $w_1,w_2,w_3,w_4,w_5$.
Observe that any pedant $S$-Steiner tree containing each vertex $w_i \ (1\leq i\leq 5)$ must occupy at
least two vertices of $\bar{S}$. Another fact is that from the
definition of $\tau_k(G)$, any pedant $S$-Steiner tree
must use at least one vertex of $\bar{S}$. So the total number of
the internally disjoint pedant $S$-Steiner trees is at
most $n-k-3$, as desired. \qed

Let $H_1=P_3\cup P_3\cup K_2$ and $H_2=P_5\cup P_3$. From Claim $1$, for $k=3$, $\bar{G}$ contains neither $H_1$ and nor $H_2$ as its subgraph. Let $H_3=P_3\cup P_3\cup K_2\cup K_1$ and $H_4=P_5\cup P_3\cup K_1$. From Claim $1$, for $k=4$, $\bar{G}$ contains neither $H_3$ and nor $H_4$ as its subgraph. Furthermore, we have the facts as follows.
\begin{itemize}
\item $\bar{G}$ contains at most two cycles.

\item $\bar{G}$ contains at most two paths of order at least $3$.

\item if $\bar{G}$ contains a cycle, then the order of this cycle is a most $7$.

\item if $\bar{G}$ contains a path of order at least $3$, then the order of this path is a most $7$.
\end{itemize}

We distinguish the following cases to show this lemma. Firstly, we
suppose that exactly two components of $\bar{G}$ are a union of two
cycles, or two paths of order at least $3$, or one is a cycle and
the other is a path of order at least $3$. Consider the case that
$\bar{G}$ contains exactly two cycles. Since $\bar{G}$ does not
contain $H_2$ as it subgraph, it follows that the order of each
cycle in $\bar{G}$ is at most $4$, and each of other components is a
isolated vertex except these cycles. Therefore,
$$
\bar{G}=C_i\cup C_j\cup (n-i-j)K_1,
$$
where $i=3,4$ and $j=3,4$. Consider the case that $\bar{G}$ contains
exactly two paths of order at least $3$. Since $\bar{G}$ does not
contain $H_2$ as it subgraph, it follows that the order of each path
in $\bar{G}$ is at most $4$, and each of other components is a
isolated vertex except these paths. Therefore,
$$
\bar{G}=P_i\cup P_j\cup (n-i-j)K_1,
$$
where $i=3,4$ and $j=3,4$. Consider the case that $\bar{G}$ contains
a cycle and a path of order at least $3$. Since $\bar{G}$ does not
contain $H_2$ as it subgraph, it follows that the order of this
cycle in $\bar{G}$ is at most $4$ and the order of this path in
$\bar{G}$ is at most $4$. Observe that each of other components is a
isolated vertex except these paths. Therefore,
$$
\bar{G}=C_i\cup P_j\cup (n-i-j)K_1,
$$
where $i=3,4$ and $j=3,4$. From the above argument, we know that $\bar{G}$ is a subgraph of $C_i\cup C_j\cup (n-i-j)K_1 \ \ (i=3,4, \ j=3,4)$, as desired.

Next, we suppose that $\bar{G}$ contains exactly one cycle or one
path of order at least $3$. Since $\bar{G}$ does not contain $H_1$
as it subgraph, if $\bar{G}$ contains exactly one cycle, then the
order of the unique cycle is at most $7$. Furthermore, if the order
of the unique cycle is $6$ or $7$, then each of other components is
a isolated vertex except this path. Therefore, $\bar{G}=C_i\cup
(n-i)K_1$, where $i=6,7$. Similarly, if the order of the unique path
is $6$ or $7$, then each of other components is a isolated vertex
except this cycle. Therefore, $\bar{G}=P_i\cup (n-i)K_1$, where
$i=6,7$. Hence, $\bar{G}$ is a subgraph of $C_i\cup (n-i)K_1 \ (i=6,7)$.

For $k=3$, we suppose the order of the unique cycle $C$ is $5$. Let
$C=v_1v_2v_3v_4v_5v_1$. If there is an independent edge $w_1w_2$ in
$\bar{G}$, then we choose $S=\{v_1,v_3,w_1\}$. Observe that any
pedant $S$-Steiner tree occupy at least one vertex in
$\bar{S}=V(G)\setminus S$. Thus there are at most $n-7$ pedant
$S$-Steiner trees by the vertices in $V(G)\setminus
\{v_1,v_2,v_3,v_4,v_5,w_1,w_2\}$. If there is no pedant $S$-Steiner
tree containing $w_2$, then any pedant $S$-Steiner tree must occupy
three vertices in $V(C)\setminus S$, and hence there are at most
$n-6$ pedant $S$-Steiner trees in $G$, a contradiction. Suppose that
there exists a pedant $S$-Steiner tree containing $w_2$, say $T$.
Then the tree $T$ must occupy at most one vertex in $V(C)\setminus
S$. Observe that there is no other pedant $S$-Steiner trees, and
hence there are at most $n-6$ pedant $S$-Steiner trees in $G$, also
a contradiction. So each of other components is a isolated vertex
except this cycle. Therefore, $\bar{G}$ is a subgraph of $C_5\cup
(n-5)K_1$. From the above argument, we know that $\bar{G}$ is a
subgraph of $C_i\cup \lfloor \frac{n-i}{2} \rfloor K_2 \ (i=3,4)$ or
$C_j\cup (n-j)K_1 \ (j=5,6,7)$ or $P_5\cup \lfloor \frac{n-5}{2}
\rfloor K_2$.

For $k=4$, we suppose the order of the unique path $P$ is $5$. Let
$P=v_1v_2v_3v_4v_5$. If there are two independent edges $w_1w_2$ and $w_3w_4$ in
$\bar{G}$, then we choose $S=\{w_1,w_3,v_2,v_4\}$. Then $|N_{\bar{G}}(S)|\geq 5$. From Claim 1,
$\tau_4(G)\leq n-7$, a contradiction. So there exists at most one nontrivial component in $\bar{G}$
except this cycle. Therefore, $\bar{G}$ is a subgraph of $C_5\cup
K_2\cup
(n-7)K_1$. From the above argument, we know that $\bar{G}$ is a
subgraph of $C_i\cup \lfloor \frac{n-i}{2} \rfloor K_2 \ (i=3,4)$ or
$C_j\cup (n-j)K_1 \ (j=6,7)$ or $C_5\cup K_2\cup (n-7)K_1$.

In the end, we suppose that $\bar{G}$ contains no cycle and no path of order at least $3$. Therefore, $\bar{G}$ is a subgraph of $\lfloor \frac{n}{2} \rfloor K_2$, and hence $\bar{G}$ is a subgraph of $C_4\cup \lfloor \frac{n-4}{2} \rfloor K_2$, as desired.

Conversely, suppose that $\bar{G}$ is a subgraph of the graphs in this lemma. For $k=3$, it suffices to show that $\tau_3(G)=n-5$. In fact, we only need to show that $\tau_3(G)=n-5$ for
\begin{itemize}
\item $\bar{G}=C_i\cup C_j\cup (n-i-j)K_1 \ \ (i=3,4, \ j=3,4)$;

\item $\bar{G}=C_i\cup \lfloor \frac{n-i}{2} \rfloor K_2  \ \ (i=3,4)$;

\item $\bar{G}=P_5\cup \lfloor \frac{n-5}{2} \rfloor K_2$;

\item $\bar{G}=C_i\cup (n-i) K_1  \ \ (i=5,6,7)$.
\end{itemize}

Suppose $\bar{G}=C_4\cup C_4\cup (n-8)K_1$. Let $C=u_1u_2u_3u_4u_1$ and $C'=v_1v_2v_3v_4v_1$ be the cycles in $\bar{G}$. It suffices to show that $\tau(S)\geq n-5$ for any $S\subseteq V(G)$ and $|S|=3$. Set $S=\{x,y,z\}$, $R=\{u_i\,|\,1\leq i\leq 4\}\cup \{v_i\,|\,1\leq i\leq 4\}$ and $V(G)\setminus R=\{w_1,w_2,\cdots,w_{n-8}\}$. If $|S\cap R|=0$, then the trees $T_i' \ (1\leq i\leq 4)$ induced by the edges in $\{xu_i,yu_i,zu_i\}$, the trees $T_i'' \ (1\leq i\leq 4)$ induced by the edges in $\{xv_i,yv_i,zv_i\}$ and the trees $T_j$ induced by the edges in $\{xw_j,yw_j,zw_j\}$ are $n-3$ pedant $S$-Steiner trees where $w_j\in \{w_1,w_2,\cdots,w_{n-8}\}\setminus \{x,y,z\}$, as desired. Suppose $|S\cap R|=1$. Without loss of generality, let $z=u_4$. Then the tree $T_2'$ induced by the edges in $\{xu_2,yu_2,zu_2\}$, the trees $T_i'' \ (1\leq i\leq 4)$ induced by the edges in $\{xv_i,yv_i,zv_i\}$ and the trees $T_j$ induced by the edges in $\{xw_j,yw_j,zw_j\}$ are $n-5$ pedant $S$-Steiner trees where $w_j\in \{w_1,w_2,\cdots,w_{n-8}\}\setminus \{x,y\}$, as desired. Suppose $|S\cap R|=2$ and $|S\cap V(C)|=2$. Without loss of generality, let $y=u_2$ and $z=u_4$. Then the trees $T_i' \ (1\leq i\leq 4)$ induced by the edges in $\{xv_i,yv_i,zv_i\}$ and the trees $T_j$ induced by the edges in $\{xw_j,yw_j,zw_j\}$ are $n-5$ pedant $S$-Steiner trees where $w_j\in \{w_1,w_2,\cdots,w_{n-8}\}\setminus \{x\}$, as desired. Suppose $|S\cap R|=2$ and $|S\cap V(C)|=|S\cap V(C')|=1$. Without loss of generality, let $y=u_4$ and $z=v_4$. Then the trees $T_1'$ induced by the edges in $\{xu_1,yv_3,v_3u_1,u_1z\}$, the trees $T_2' \ (1\leq i\leq 4)$ induced by the edges in $\{xv_1,yv_1,u_3v_1,u_3z\}$ and the trees $T_j$ induced by the edges in $\{xw_j,yw_j,zw_j\}$ are $n-5$ pedant $S$-Steiner trees where $w_j\in \{w_1,w_2,\cdots,w_{n-8}\}\setminus \{x\}$, as desired.
Suppose $|S\cap R|=3$ and $|S\cap V(C)|=3$. Then the trees $T_i' \ (1\leq i\leq 4)$ induced by the edges in $\{xv_i,yv_i,zv_i\}$ and the trees $T_j$ induced by the edges in $\{xw_j,yw_j,zw_j\}$ are $n-4$ pedant $S$-Steiner trees where $w_j\in \{w_1,w_2,\cdots,w_{n-8}\}$, as desired. Suppose $|S\cap R|=3$ and $|S\cap V(C)|=2$. Without loss of generality, let $x=u_2$, $y=u_4$ and $z=v_4$. Then the tree $T_1'$ induced by the edges in $\{xv_3,yv_3,v_3u_3,u_3z\}$, the tree $T_2'$ induced by the edges in $\{xv_1,yv_1,u_1v_1,u_1z\}$ and the trees $T_j$ induced by the edges in $\{xw_j,yw_j,zw_j\}$ are $n-4$ pedant $S$-Steiner trees where $w_j\in \{w_1,w_2,\cdots,w_{n-7}\}$, as desired.

Suppose $\bar{G}=C_7\cup (n-7)K_1$. Let $C=u_1u_2\cdots u_7u_1$ be
the cycle in $\bar{G}$. It suffices to show that $\tau(S)\geq n-5$
for any $S\subseteq V(G)$ and $|S|=3$. Set $S=\{x,y,z\}$ and
$V(G)\setminus V(C)=\{w_1,w_2,\cdots,w_{n-7}\}$. If $|S\cap
V(C)|=0$, then the trees $T_i' \ (1\leq i\leq 7)$ induced by the
edges in $\{xu_i,yu_i,zu_i\}$ and the trees $T_j$ induced by the
edges in $\{xw_j,yw_j,zw_j\}$ are $n-3$ pedant $S$-Steiner trees
where $w_j\in \{w_1,w_2,\cdots,w_{n-7}\}\setminus \{x,y,z\}$, as
desired. Suppose $|S\cap V(C)|=1$. Without loss of generality, let
$z=u_7$. Then the trees $T_i' \ (2\leq i\leq 5)$ induced by the
edges in $\{xu_i,yu_i,zu_i\}$ and the trees $T_j$ induced by the
edges in $\{xw_j,yw_j,zw_j\}$ are $n-5$ pedant $S$-Steiner trees
where $w_j\in \{w_1,w_2,\cdots,w_{n-7}\}\setminus \{x,y\}$, as
desired. Suppose that $|S\cap V(C)|=3$. We only need to consider the
cases $\{x,y,z\}=\{u_1,u_2,u_3\}$, $\{x,y,z\}=\{u_1,u_2,u_4\}$,
$\{x,y,z\}=\{u_1,u_3,u_5\}$ and $\{x,y,z\}=\{u_1,u_2,u_5\}$. If
$\{x,y,z\}=\{u_1,u_2,u_3\}$, then the trees $T_i' \ (i=5,6)$ induced
by the edges in $\{xu_i,yu_i,zu_i\}$ and the trees $T_j$ induced by
the edges in $\{xw_j,yw_j,zw_j\}$ are $n-5$ pedant $S$-Steiner trees
where $w_j\in \{w_1,w_2,\cdots,w_{n-7}\}$, as desired. If
$\{x,y,z\}=\{u_1,u_2,u_4\}$, then the tree $T_1'$ induced by the edges in
$\{xu_6,yu_6,zu_6\}$, the tree $T_2'$ induced by the edges in
$\{u_2u_7,u_4u_7,u_5u_7,u_1u_5\}$ and the trees $T_j$ induced by the
edges in $\{xw_j,yw_j,zw_j\}$ are $n-5$ pedant $S$-Steiner trees
where $w_j\in \{w_1,w_2,\cdots,w_{n-7}\}$, as desired. If
$\{x,y,z\}=\{u_1,u_3,u_5\}$, then the tree induced by the edges in
$\{u_1u_4,u_4u_7,u_3u_7,u_7u_5\}$, the tree induced by the edges in
$\{u_1u_6,u_2u_6,u_3u_6,u_2u_5\}$ and the trees $T_j$ induced by the
edges in $\{xw_j,yw_j,zw_j\}$ are $n-5$ pedant $S$-Steiner trees
where $w_j\in \{w_1,w_2,\cdots,w_{n-7}\}$, as desired. If
$\{x,y,z\}=\{u_1,u_2,u_5\}$, then the trees induced by the edges in
$\{u_1u_6,u_2u_6,u_3u_6,u_3u_5\}$, the tree induced by the edges in
$\{u_2u_4,u_1u_4,u_4u_7,u_7u_5\}$ and the trees $T_j$ induced by the
edges in $\{xw_j,yw_j,zw_j\}$ are $n-5$ pedant $S$-Steiner trees
where $w_j\in \{w_1,w_2,\cdots,w_{n-7}\}$, as desired. Suppose
$|S\cap V(C)|=2$. We only need to consider the cases
$\{x,y,z\}=\{w_1,u_3,u_4\}$, $\{x,y,z\}=\{w_1,u_2,u_4\}$ and
$\{x,y,z\}=\{w_1,u_2,u_5\}$. If $\{x,y,z\}=\{w_1,u_2,u_3\}$, then
the trees $T_i' \ (i=1,6,7)$ induced by the edges in
$\{xu_i,yu_i,zu_i\}$ and the trees $T_j$ induced by the edges in
$\{xw_j,yw_j,zw_j\}$ are $n-5$ pedant $S$-Steiner trees where
$w_j\in \{w_2,w_3,\cdots,w_{n-7}\}$, as desired. If
$\{x,y,z\}=\{w_1,u_2,u_4\}$, then the trees $T_i' \ (i=6,7)$ induced
by the edges in $\{xu_i,yu_i,zu_i\}$, the tree $T''$ induced by the
edges in $\{xu_1,u_1u_5,yu_4,zu_1\}$ and the trees $T_j$ induced by
the edges in $\{xw_j,yw_j,zw_j\}$ are $n-5$ pedant $S$-Steiner trees
where $w_j\in \{w_2,w_3,\cdots,w_{n-7}\}$, as desired. If
$\{x,y,z\}=\{w_1,u_2,u_5\}$, then the tree $T_1'$ induced by the
edges in $\{xu_7,yu_7,zu_7\}$, the tree $T_2'$ induced by the edges
in $\{w_1u_4,u_1u_4,u_2u_4,u_1u_5\}$, the tree $T_3'$ induced by the edges
in $\{w_1u_3,u_3u_5,u_3u_6,u_2u_6\}$ and the trees $T_j$ induced by the
edges in $\{xw_j,yw_j,zw_j\}$ are $n-5$ pedant $S$-Steiner trees
where $w_j\in \{w_2,w_3,\cdots,w_{n-7}\}$, as desired.

Suppose $\bar{G}=P_5\cup \lfloor \frac{n-5}{2}\rfloor K_2$. Let $P=u_1u_2u_3u_4u_5$ be the unique path in $\bar{G}$. It suffices to show that $\tau(S)\geq n-5$ for any $S\subseteq V(G)$ and $|S|=3$. Set $S=\{x,y,z\}$ and $V(G)\setminus V(P)=\{w_1,w_2,\cdots,w_{n-5}\}$. Note that $C_7\cup \lfloor \frac{n-7}{2}\rfloor K_2$ is a subgraph of $\bar{G}$. Suppose that $|S\cap V(P)|=3$ or $|S\cap V(P)|=2$. Recall that we have checked that there are $n-5$ pedant $S$-Steiner trees in the complement of $C_7\cup (n-7) K_1$. In fact, one can check that if $|S\cap V(P)|=3$ or $|S\cap V(P)|=2$, then there are $n-5$ pedant $S$-Steiner trees in the complement of $C_7\cup \lfloor \frac{n-7}{2}\rfloor K_2$, and hence there are $n-5$ pedant $S$-Steiner trees in the complement of $\bar{G}$. We now assume that $|S\cap V(P_5)|=1$ or $|S\cap V(P)|=0$. If $|S\cap V(P)|=0$, then the trees $T_i' \ (1\leq i\leq 5)$ induced by the edges in $\{xu_i,yu_i,zu_i\}$ and the trees $T_j$ induced by the edges in $\{xw_j,yw_j,zw_j\}$ are $n-4$ pedant $S$-Steiner trees where $w_j\in \{w_1,w_2,\cdots,w_{n-5}\}\setminus \{x,y,z,w_1\}$ where $w_1$ is a vertex adjacent to $x$ in $\bar{G}$, as desired. Suppose $|S\cap V(P)|=1$. Without loss of generality, let $z\in V(P)$ and $x,y\in V(G)\setminus V(P)$. If $xy\in E(\bar{G})$, then one can check that there are $n-5$ pedant $S$-Steiner trees in $G$ since $C_7\cup \lfloor \frac{n-7}{2}\rfloor K_2$ is a subgraph of $\bar{G}$. We may assume $xy\notin E(\bar{G})$. Without loss of generality, let $x=w_1$, $y=w_3$, $w_1w_2\in E(\bar{G})$ and $w_3w_4\in E(\bar{G})$. If $z=u_1$, then the trees $T_i' \ (3\leq i\leq 5)$ induced by the edges in $\{xu_i,yu_i,zu_i\}$, the tree induced by the edges in $\{u_1w_2,w_2w_4,w_1w_4,w_2w_3\}$ and the trees $T_j \ (5\leq i\leq n-5)$ induced by the edges in $\{xw_j,yw_j,zw_j\}$ are $n-5$ pedant $S$-Steiner trees, as desired.
If $z=u_2$, then the trees $T_i' \ (4\leq i\leq 5)$ induced by the edges in $\{xu_i,yu_i,zu_i\}$, the tree induced by the edges in $\{u_2w_4,v_3w_4,u_3w_3,v_3w_1\}$, the tree induced by the edges in $\{u_2w_2,u_1w_2,w_2w_3,u_1w_1\}$ and the trees $T_j \ (5\leq i\leq n-5)$ induced by the edges in $\{xw_j,yw_j,zw_j\}$ are $n-5$ pedant $S$-Steiner trees, as desired. If $z=u_3$, then the trees $T_i' \ (i=1,5)$ induced by the edges in $\{xu_i,yu_i,zu_i\}$, the tree induced by the edges in $\{u_4w_3,w_1u_4,u_4w_4,w_4u_3\}$, the tree induced by the edges in $\{u_2w_1,w_2u_2,u_3w_2,w_2w_3\}$ and the trees $T_j$ induced by the edges in $\{xw_j,yw_j,zw_j\}$ are $n-5$ pedant $S$-Steiner trees, as desired.

For other cases, one can also check that $\tau(S)\geq n-5$ for any $S\subseteq V(G)$ and $|S|=3$. Therefore, $\tau_3(G)\geq n-5$. From Theorem \ref{th5}, we have $\tau_3(G)=n-5$.

For $k=4$, it suffices to show that $\tau_4(G)=n-6$. In fact, we only need to show that $\tau_4(G)=n-6$ for
\begin{itemize}
\item $\bar{G}=C_i\cup C_j\cup (n-i-j)K_1 \ \ (i=3,4, \ j=3,4)$;

\item $\bar{G}=C_i\cup \lfloor \frac{n-i}{2} \rfloor K_2  \ \ (i=3,4)$;

\item $\bar{G}=P_5\cup K_2\cup (n-7) K_1$;

\item $\bar{G}=C_i\cup (n-i) K_1  \ \ (i=5,6,7)$.
\end{itemize}

Suppose $\bar{G}=C_4\cup C_4\cup (n-8)K_1$. Let $C=u_1u_2u_3u_4u_1$ and $C'=v_1v_2v_3v_4v_1$ be the cycles in $\bar{G}$. It suffices to show that $\tau(S)\geq n-6$ for any $S\subseteq V(G)$ and $|S|=4$. Set $S=\{x,y,z,r\}$, $R=\{u_i\,|\,1\leq i\leq 4\}\cup \{v_i\,|\,1\leq i\leq 4\}$ and $V(G)\setminus R=\{w_1,w_2,\cdots,w_{n-8}\}$.
If $|S\cap R|=0$, then the trees $T_i' \ (1\leq i\leq 4)$ induced by the edges in $\{xu_i,yu_i,zu_i,ru_i\}$, the trees $T_i'' \ (1\leq i\leq 4)$ induced by the edges in $\{xv_i,yv_i,zv_i,rv_i\}$ and the trees $T_j$ induced by the edges in $\{xw_j,yw_j,zw_j,rw_j\}$ are $n-4$ pedant $S$-Steiner trees where $w_j\in \{w_1,w_2,\cdots,w_{n-8}\}\setminus \{x,y,z,r\}$, as desired.
Suppose $|S\cap R|=1$. Without loss of generality, let $r=v_3$. Then the tree $T'$ induced by the edges in $\{xv_1,yv_1,zv_1,rv_1\}$, the trees $T_i'' \ (1\leq i\leq 4)$ induced by the edges in $\{xu_i,yu_i,zu_i,zu_i\}$ and the trees $T_j$ induced by the edges in $\{xw_j,yw_j,zw_j,rw_j\}$ are $n-6$ pedant $S$-Steiner trees where $w_j\in \{w_1,w_2,\cdots,w_{n-8}\}\setminus \{x,y,z\}$, as desired.
Suppose $|S\cap R|=2$ and $|S\cap V(C')|=2$. Without loss of generality, let $r=v_1$ and $z=v_3$. Then the trees $T_i' \ (1\leq i\leq 4)$ induced by the edges in $\{xu_i,yu_i,zu_i,ru_i\}$ and the trees $T_j$ induced by the edges in $\{xw_j,yw_j,zw_j,rw_j\}$ are $n-6$ pedant $S$-Steiner trees where $w_j\in \{w_1,w_2,\cdots,w_{n-8}\}\setminus \{x,y\}$, as desired.
Suppose $|S\cap R|=2$ and $|S\cap V(C)|=|S\cap V(C')|=1$. Without loss of generality, let $r=u_1$ and $z=v_3$. Then the tree $T_1'$ induced by the edges in $\{xu_2,yu_2,v_2u_2,u_2z,v_2r\}$, the tree $T_2'$ induced by the edges in $\{xv_4,yv_4,rv_4,u_4v_4,u_4z\}$, the tree $T_3'$ induced by the edges in $\{xv_1,yv_1,zv_1,rv_1\}$, the tree $T_4'$ induced by the edges in $\{xu_3,yu_3,zu_3,ru_3\}$ and the trees $T_j$ induced by the edges in $\{xw_j,yw_j,zw_j,rw_j\}$ are $n-6$ pedant $S$-Steiner trees where $w_j\in \{w_1,w_2,\cdots,w_{n-8}\}\setminus \{x,y\}$, as desired.
Suppose $|S\cap R|=3$ and $|S\cap V(C)|=3$. Without loss of generality, let $S\cap V(C)=\{y,z,r\}$. Then the trees $T_i' \ (1\leq i\leq 4)$ induced by the edges in $\{xv_i,yv_i,zv_i,rv_i\}$ and the trees $T_j$ induced by the edges in $\{xw_j,yw_j,zw_j,rw_j\}$ are $n-6$ pedant $S$-Steiner trees where $w_j\in \{w_1,w_2,\cdots,w_{n-8}\}\setminus \{x\}$, as desired.
Suppose $|S\cap R|=3$ and $|S\cap V(C)|=2$. Without loss of generality, let $r=v_4$, $y=u_3$ and $z=u_4$.
Then the trees $T_1'$ induced by the edges in $\{xv_3,yv_3,zv_3,v_3u_2,u_2r\}$, the trees $T_2'$ induced by the edges in $\{xu_1,yu_1,u_1r,u_1v_1,v_1z\}$ and the trees $T_j$ induced by the edges in $\{xw_j,yw_j,zw_j,rw_j\}$ are $n-6$ pedant $S$-Steiner trees where $w_j\in \{w_1,w_2,\cdots,w_{n-8}\}\setminus \{x\}$, as desired.
Suppose $|S\cap R|=4$ and $|S\cap V(C)|=4$. Without loss of generality, let $S\cap V(C)=\{x,y,z,r\}$. Then the trees $T_i' \ (1\leq i\leq 4)$ induced by the edges in $\{xv_i,yv_i,zv_i,rv_i\}$ and the trees $T_j$ induced by the edges in $\{xw_j,yw_j,zw_j,rw_j\}$ are $n-6$ pedant $S$-Steiner trees where $w_j\in \{w_1,w_2,\cdots,w_{n-8}\}$, as desired.
Suppose $|S\cap R|=4$ and $|S\cap V(C)|=3$. Without loss of generality, let $S\cap V(C)=\{x,y,z\}$. Then the tree $T_1'$ induced by the edges in $\{xv_2,yv_2,zv_2,rv_2\}$, the tree $T_2'$ induced by the edges in $\{xv_1,yv_1,zv_1,u_1v_1,ru_1\}$
and the trees $T_j$ induced by the edges in $\{xw_j,yw_j,zw_j,rw_j\}$ are $n-6$ pedant $S$-Steiner trees where $w_j\in \{w_1,w_2,\cdots,w_{n-8}\}$, as desired.
Suppose $|S\cap R|=4$ and $|S\cap V(C)|=2$. Without loss of generality, let $x=u_2$, $y=u_4$, $z=v_2$ and $r=v_4$.
Then the trees $T_1'$ induced by the edges in $\{xv_1,yv_1,v_1u_3,zu_3,u_3r\}$, the trees $T_2'$ induced by the edges in $\{xv_3,yv_3,u_1v_3,u_1r,u_1z\}$ and the trees $T_j$ induced by the edges in $\{xw_j,yw_j,zw_j,rw_j\}$ are $n-6$ pedant $S$-Steiner trees where $w_j\in \{w_1,w_2,\cdots,w_{n-8}\}$, as desired.

Suppose $\bar{G}=C_7\cup (n-7)K_1$. Let $C=u_1u_2\cdots u_7u_1$ be
the cycle in $\bar{G}$. It suffices to show that $\tau(S)\geq n-6$
for any $S\subseteq V(G)$ and $|S|=4$. Set $S=\{x,y,z,r\}$ and
$V(G)\setminus V(C)=\{w_1,w_2,\cdots,w_{n-7}\}$.
If $|S\cap
V(C)|=0$, then the trees $T_i' \ (1\leq i\leq 7)$ induced by the
edges in $\{xu_i,yu_i,zu_i,ru_i\}$ and the trees $T_j$ induced by the
edges in $\{xw_j,yw_j,zw_j,rw_j\}$ are $n-4$ pedant $S$-Steiner trees
where $w_j\in \{w_1,w_2,\cdots,w_{n-7}\}\setminus \{x,y,z,r\}$, as
desired. Suppose $|S\cap V(C)|=1$. Without loss of generality, let
$z=u_1$. Then the trees $T_i' \ (3\leq i\leq 6)$ induced by the
edges in $\{xu_i,yu_i,zu_i,ru_i\}$ and the trees $T_j$ induced by the
edges in $\{xw_j,yw_j,zw_j,rw_j\}$ are $n-6$ pedant $S$-Steiner trees
where $w_j\in \{w_1,w_2,\cdots,w_{n-7}\}\setminus \{x,y,z\}$, as
desired. Suppose $|S\cap V(C)|=2$. Without loss of generality, let $z=u_1$ and $r=u_4$. Then the tree $T_1'$ induced by the edges in
$\{xu_6,yu_6,zu_6,ru_6\}$, the tree $T_2'$ induced by the edges in
$\{xu_3,yu_3,zu_3,u_3u_7,ru_7\}$, the tree $T_3'$ induced by the edges in
$\{xu_2,yu_2,ru_2,u_2u_5,zu_5\}$ and the trees $T_j$ induced by the edges in
$\{xw_j,yw_j,zw_j,rw_j\}$ are $n-6$ pedant $S$-Steiner trees where
$w_j\in \{w_2,w_3,\cdots,w_{n-7}\}\setminus \{x,y\}$, as desired.
Suppose that $|S\cap V(C)|=3$. Without loss of generality, let $y=u_1$, $z=u_4$ and $r=u_6$. Then the tree $T'$ induced by the edges in
$\{xu_5,yu_5,u_5u_2,zu_2,ru_2\}$, the tree $T''$ induced by the edges in
$\{xu_3,yu_3,zu_7,u_3u_7,ru_3\}$ and the trees $T_j$ induced by the edges in
$\{xw_j,yw_j,zw_j,rw_j\}$ are $n-6$ pedant $S$-Steiner trees where
$w_j\in \{w_2,w_3,\cdots,w_{n-7}\}\setminus \{x\}$, as desired.
Suppose that $|S\cap V(C)|=4$. Without loss of generality, let $x=u_1$, $y=u_3$, $z=u_5$ and $r=u_7$.
Then the tree $T'$ induced by the edges in
$\{xu_6,yu_6,u_6u_2,zu_2,ru_2\}$ and the trees $T_j$ induced by the edges in
$\{xw_j,yw_j,zw_j,rw_j\}$ are $n-6$ pedant $S$-Steiner trees where
$w_j\in \{w_2,w_3,\cdots,w_{n-7}\}$, as desired.

For other cases, one can also check that $\tau(S)\geq n-6$ for any $S\subseteq V(G)$ and $|S|=4$. Therefore, $\tau_4(G)\geq n-6$. From Theorem \ref{th5}, we have $\tau_4(G)=n-6$.\qed
\end{pf}\\

For $3\leq k\leq n$, graphs with $\tau_k(G)=n-k-2$ can also be characterized in the following.

\begin{thm}\label{th6}
Let $k,n$ be two integers with $3\leq k\leq n$ and $n\geq 15$, and let $G$ be a connected graph.

$(1)$ For $5\leq k\leq n$, $\tau_k(G)=n-k-2$ if and only if $\bar{G}=P_3\cup (n-3)K_1$ or
$G$ satisfies all the following conditions.

$\bullet$ $1\leq \Delta(\bar{G})\leq 2$;

$\bullet$ $e(\bar{G})\geq 3$;

$\bullet$ for any $R\subseteq V(G)$ with $|E_{\bar{G}}[w,\bar{R}]|\geq 1$ where $w\in R$ and $\bar{R}=V(G)\setminus R$, the size of $R$ is at most $4$.

$(2)$ For $k=3$, $\tau_k(G)=n-k-2$ if and only if $\bar{G}$ is a subgraph of one of the following graphs.
\begin{itemize}
\item $C_i\cup C_j\cup (n-i-j)K_1 \ \ (i=3,4, \ j=3,4)$;

\item $C_i\cup \lfloor \frac{n-i}{2} \rfloor K_2  \ \ (i=3,4)$;

\item $P_5\cup \lfloor \frac{n-5}{2} \rfloor K_2$;

\item $C_i\cup (n-i) K_1  \ \ (i=5,6,7)$.
\end{itemize}

$(3)$ For $k=4$, $\tau_k(G)=n-k-2$ if and only if $\bar{G}$ is a subgraph of one of the following graphs.
\begin{itemize}
\item $C_i\cup C_j\cup (n-i-j)K_1 \ \ (i=3,4, \ j=3,4)$;

\item $C_i\cup \lfloor \frac{n-i}{2} \rfloor K_2  \ \ (i=3,4)$;

\item $C_5\cup K_2\cup (n-7) K_1$;

\item $C_i\cup (n-i) K_1  \ \ (i=6,7)$.
\end{itemize}
\end{thm}
\begin{pf}
From Lemma \ref{lem6}, the conclusion is true for the case $k=3,4$. We now assume $5\leq k\leq n$. Suppose $\tau_k(G)=n-k-2$.
If $e(\bar{G})=2$, then $\bar{G}=P_3\cup (n-3)K_1$ by Theorem \ref{th5}. Conversely, we suppose $\bar{G}=P_3\cup (n-3)K_1$. One can check that there exist $n-k-2$ pedant $S$-Steiner tree in $G$ for any $S\subseteq V(G)$ and $|S|=k$. Then
$\tau_k(G)\geq n-k-2$. From Theorem \ref{th5}, we have $\tau_k(G)=n-k-2$.

From now on, we assume $e(\bar{G})\geq 3$. Suppose $\tau_k(G)=n-k-2$. From Lemma \ref{lem1}, we have
$\delta(G)\geq k+(n-k-2)-1=n-3$. Therefore,
$\Delta(\bar{G})=n-1-\delta(G)\leq 2$. Combining this with Theorems
\ref{th4} and \ref{th5}, we have $1\leq \Delta(\bar{G})\leq 2$.
 Since $e(\bar{G})\geq 3$, it follows
that each component of $\bar{G}$ is a path or a cycle (Note that an
isolated vertex can be seen a path of order $1$). From Lemma \ref{lem6}, the result is true for $k=3,4$. For $5\leq k\leq n$, we
have the following claim.

\textbf{Claim 1.} For any $R\subseteq V(G)$ with $|E_{\bar{G}}[w,\bar{R}]|\geq 1$ where $w\in R$ and $\bar{R}=V(G)\setminus R$, $|R|\leq 4$.

\noindent{\bf Proof of Claim 1.} Assume, to the contrary, that $|R|=5$. Set
$R=\{w_i\,|\,1\leq i\leq 5\}$, and let $U$ be the vertex set such that $|E_{\bar{G}}[w_i,U]|\geq 1$ and $|U|\leq 5$.
Choose $S\subseteq V(G)\setminus R$ and $|S|=k$ such that $S$ contains the vertex set $U$.
Let $S=\{u_1,u_2,\cdots,u_k\}$ and $\bar{S}=V(G)\setminus S=\{w_1,w_2,\cdots,w_{n-k}\}$.
Note that $|E_{\bar{G}}[w_i,\bar{R}]|\geq 1$ for any $w_i \ (1\leq i\leq 5)$.
Observe that any pedant $S$-Steiner tree containing each vertex $w_i \ (1\leq i\leq 5)$ must occupy at
least two vertices of $\bar{S}$. Another fact is that from the
definition of $\tau_k(G)$, any pedant $S$-Steiner tree
must use at least one vertex of $\bar{S}$. So the total number of
the internally disjoint pedant $S$-Steiner trees is at
most $2+(n-k-5)=n-k-3$, a contradiction. So $|R|\leq 4$ for
$5\leq k\leq n$.

Conversely, suppose $G$ satisfies the condition of this theorem for $5\leq k\leq n$. It is
clear that we only need to prove that $\tau_k(G)\geq n-k-2$ where
$G=K_n\setminus M$ such that $1\leq \Delta(\bar{G})\leq 2$ and for any $R\subseteq V(G)$ with $|E_{\bar{G}}[w,\bar{R}]|\geq 1$ where $w\in R$ and $\bar{R}=V(G)\setminus R$, the size of $R$ is exactly $4$. Set $R=\{w_i\,|\,1\leq i\leq 4\}$.
From the definition of $\tau_k(G)$, it suffices to show that
$\tau(S)\geq n-k-2$ for any $S\subseteq V(G)$ with $|S|=k$.
Since $|R|=4$, we can assume $w_ju_i\in E(G)$ for any $1\leq i\leq k$ and $5\leq j\leq n-k$. Then the
trees $T_j$ induced by the edges in $\{u_ju_1,u_ju_2,\cdots,u_ju_k\} \ (5\leq j\leq
n-k)$ form $n-k-4$ pedant $S$-Steiner trees.
Recall that for any $R\subseteq V(G)$ with $|E_{\bar{G}}[w,\bar{R}]|\geq 1$ where $w\in R$ and $\bar{R}=V(G)\setminus R$, $|R|\leq 4$. Therefore, there are at most four vertices in $S$, without loss of generality, say $u_1,u_2,u_3,u_4$ such that $|E_{\bar{G}}[u_i,\bar{R}]|\geq 1$ for any $u_i \ (1\leq i\leq
4)$. Then $w_ju_i\in E(G)$ for any $4\leq i\leq k$ and $1\leq j\leq 4$.
Since $1\leq \Delta(\bar{G})\leq 2$, it follows that we may assume that $d_{\bar{G}}(w_j)=2$ for $1\leq j\leq 4$. Without loss of generality, let $w_1u_1,w_1u_2\in M$. Since $\Delta(\bar{G})\leq 2$, it follows that there exists a vertex in $\{w_2,w_3,w_4\}$, say $w_4$, such that $w_4u_1\notin M$ and $w_4u_2\notin M$. Furthermore, $w_4u_1\in E(G)$, $w_4u_2\in E(G)$ and hence the
trees $T_1$ induced by the edges in $\{w_4u_1,w_4u_2,w_4w_1,w_1u_3\cdots,w_1u_k\}$ is a pedant $S$-Steiner tree. Since $\Delta(\bar{G})\leq 2$, we only need to consider the case $u_1w_2,u_2w_3\in M$ and the case $u_1w_2,u_2w_2\in M$. For the former case, the
trees $T_2$ induced by the edges in $\{w_3u_1,w_2w_3,w_2u_2\cdots,w_2u_k\}$ is a pedant $S$-Steiner tree. For the latter case, the
trees $T_2$ induced by the edges in $\{w_3u_1,w_3u_2,w_2w_3,w_3u_3\cdots,w_2u_k\}$ is a pedant $S$-Steiner tree. Therefore, the trees $T_1,T_2$ together with the trees $T_5,\cdots,T_{n-k}$ are $n-k-2$ internally disjoint pedant $S$-Steiner trees.
From the above argument, we conclude that for any $S\subseteq V(G)$
with $|S|=k$, $\tau(S)\geq n-k-2$. From the arbitrariness of $S$, we
have $\tau_k(G)\geq n-k-2$. Combining this with Theorem
\ref{th5}, $\tau_k(G)=n-k-2$. \qed
\end{pf}\\

If $k=n-1$, then $0\leq \tau_{n-1}(G)\leq 1$ by Proposition
\ref{pro1}.

\begin{cor}\label{cor2}
Let $G$ be a connected graph of order $n$. Then

$(1)$ $\tau_{n-1}(G)=1$ if and only if $G$ is a complete graph of
order $n$.

$(2)$ $\tau_{n-1}(G)=0$ if and only if $G$ is not a complete graph
of order $n$.
\end{cor}

If $k=n-2$, then $0\leq \tau_{n-2}(G)\leq 2$ by Proposition
\ref{pro1}.

\begin{cor}\label{cor3}
Let $G$ be a connected graph of order $n$. Then

$(1)$ $\tau_{n-2}(G)=2$ if and only if $G$ is a complete graph of
order $n$.

$(2)$ $\tau_{n-2}(G)=1$ if and only if $G=K_n\setminus M$ and $1\leq
|M|\leq 2$, where $M$ is a matching of $K_n$ for $n\geq 7$.

$(3)$ $\tau_{n-2}(G)=0$ if and only if $G$ is one of the other
graphs.
\end{cor}

\section{Graphs with small pedant tree-connectivity}

Given a vertex $x$ and a set $U$ of vertices, an \emph{$(x,U)$-fan}
is a set of paths from $x$ to $U$ such that any two of them share
only the vertex $x$. The size of an $(x,U)$-fan is the number of
internally disjoint paths from $x$ to $U$.

\begin{lem}{\upshape (Fan Lemma, \cite{West}, p-170)}\label{lem7}
A graph is $k$-connected if and only if it has at least $k+1$
vertices and, for every choice of $x$, $U$ with $|U|\geq k$, it has
an $(x,U)$-fan of size $k$.
\end{lem}

We now turn our attention to characterize graphs with
$\tau_k(G)=0$.

\begin{thm}\label{th7}
Let $k,n$ be two integers with $3\leq k\leq n$. Let $G$ be a
connected graph of order $n$. Then $\tau_k(G)=0$ if and only if $G$
satisfies one of the following conditions.

$(1)$ $0\leq \kappa(G)\leq k-2$;

$(2)$ $\kappa(G)=\delta(G)=k-1$;

$(3)$ $\kappa(G)=k-1$, $\delta(G)\geq k$, and there exists a vertex
subset $S$ of $V(G)$ with $|S|=k$ such that
for any $S'\subseteq S$ with $|S'|=k-1$,
for any vertex $x\in V(G_i)\setminus S$, and any
$(x,S')$-fan, $u_1$ is an internal vertex of some path belonging to
this $(x,S')$-fan, where $G_i$ is the connected component of $G\setminus S'$ containing $u_1$.
\end{thm}
\begin{pf}
Suppose $\tau_k(G)=0$. If $k=n$, then $\tau_n(G)=0$
if and only if $G$ is a connected graph by Corollary \ref{cor1},
which implies that $\tau_n(G)=0$ if and only if $0\leq \kappa(G)\leq
n-2$ or $\kappa(G)=\delta(G)=n-1$, as desired. We now assume $3\leq
k\leq n-1$. Then we have the following claim.

{\bf Claim 1.} $\kappa(G)\leq k-1$.

\noindent{\bf Proof of Claim 1.} Assume, to the contrary, that $\kappa(G)\geq k$. For any $S\subseteq V(G)$ and
$|S|=k$, there is a vertex $x\in V(G)\setminus S$ since $3\leq k\leq
n-1$. Let $S=\{u_1,u_2,\cdots,u_k\}$. Since $\kappa(G)\geq k$, it
follows from Lemma \ref{lem7} that there exists an $(x,S)$-fan of size $k$ in $G$. Let
$P_1,P_2,\cdots,P_k$ be the $k$ internally disjoint paths of this
$(x,S)$-fan. Then the tree $T$ induced by the edges in $E(P_1)\cup E(P_2)\cup \cdots \cup E(P_k)$ is a
pedant $S$-Steiner tree, which implies that $\tau(S)\geq
1$. From the arbitrariness of $S$, we have $\tau_k(G)\geq 1$, a
contradiction.\qed

From Claim 1, we have $0\leq \kappa(G)\leq k-1$. If $0\leq
\kappa(G)\leq k-2$, then $(1)$ holds. If $\kappa(G)=k-1$, then
$\delta(G)\geq \kappa(G)=k-1$. Furthermore, if
$\kappa(G)=\delta(G)=k-1$, then $(2)$ holds.
From now on, we assume that $\kappa(G)=k-1$ and $\delta(G)\geq k$.

{\bf Claim 2.} There exists a vertex
subset $S$ of $V(G)$ with $|S|=k$ such that for any $S'\subseteq S$ with $|S'|=k-1$, if $S'$ is not a vertex
cut set of $G$, then for any vertex $x\in V(G)\setminus S$, and any
$(x,S')$-fan, $u_1$ is an internal vertex of a path belonging to
this $(x,S')$-fan.

\noindent{\bf Proof of Claim 2.} Assume, to the contrary, that
for any $S\subseteq V(G)$ with $|S|=k$, there exists a vertex
subset $S'$ in $S$ with $|S'|=k-1$ such that $S'$ is not a vertex
cut of $G$, and there exists a vertex $x\in V(G)\setminus S$ and an
$(x,S')$-fan, $u_1$ is not an internal vertex of for any path of
this $(x,S')$-fan.
Let $S'=\{u_2,u_3,\cdots,u_k\}$ and $S\setminus
S'=\{u_1\}$, where $S=\{u_1,u_2,\cdots,u_k\}$.
Since
$\kappa(G)=k-1$, it follows from Lemma \ref{lem7} that there is an $(x,S')$-fan in $G$,
where $x\in V(G)\setminus S$. Note that $x\neq u_1$. Denote by $P_2,P_3,\cdots,P_{k-1}$ the
$k$ internally disjoint paths connecting $x$ and
$u_2,u_3,\cdots,u_{k-1}$ of this $(x,S')$-fan, respectively. Recall that
$u_1$ is not an internal vertex of for any path of
this $(x,S')$-fan. Since
$\kappa(G)=k-1$ and $S'$ is not a vertex cut of $G$, it follows that
$G\setminus S'$ is connected and hence there is a path connecting
$x$ and $u_1$, say $P_1$. Clearly, the graph $H$ induced by the edges in $E(P_1)\cup E(P_2)\cup \cdots \cup E(P_k)$ contains a pedant $S$-Steiner tree,
which implies that $\tau(S)\geq 1$. From the arbitrariness of $S$,
we have $\tau_k(G)\geq 1$, a contradiction.\qed

Furthermore, we have the following claim.

{\bf Claim 3.} There exists a vertex
subset $S$ of $V(G)$ with $|S|=k$ such that for any $S'\subseteq S$ with $|S'|=k-1$, if $S'$ is a vertex
cut set of $G$, then for any vertex $x\in V(G_i)\setminus u_1$, and any
$(x,S')$-fan, $u_1$ is an internal vertex of a path belonging to
this $(x,S')$-fan, where $G_i$ is the connected component of $G\setminus S'$ containing $u_1$.

\noindent{\bf Proof of Claim 3.} Assume, to the contrary, that for any $S\subseteq V(G)$ with $|S|=k$, there exists a vertex
subset $S'$ in $S$ with $|S'|=k-1$ such that $S'$ is a vertex
cut of $G$, and
there exist a vertex $x\in V(G_i)\setminus u_1$ and an
$(x,S')$-fan such that $u_1$ does not belong to this $(x,S')$-fan,
where $G_i$ is the connected component of $G\setminus S'$ containing $u_1$.
Denote by $P_2,P_3,\cdots,P_{k-1}$ the $k$ internally disjoint paths
connecting $x$ and $u_2,u_3,\cdots,u_{k-1}$ of this $(x,S')$-fan,
respectively.
Since $G_i$ is connected, there is
a path connecting $x$ and $u_1$, say $P_1$. Clearly, the graph
$H$ induced by the edges in $E(P_1)\cup E(P_2)\cup \cdots \cup E(P_k)$ contains a pedant $S$-Steiner tree, which implies that $\tau(S)\geq 1$.
From the
arbitrariness of $S$, we have $\tau_k(G)\geq 1$, a contradiction.\qed

From Claims 2 and 3, we know that
$(3)$ holds.

Conversely, we suppose that $G$ is a graph satisfying one of conditions in this
theorem. Our aim is to show $\tau_k(G)=0$. Suppose $0\leq \kappa(G)\leq k-2$. If $\tau_k(G)\geq 1$,
then we have $\kappa(G)\geq k-1$ by Lemma \ref{lem1}, a
contradiction. Therefore, $\tau_k(G)=0$, as desired. Suppose
$\kappa(G)=\delta(G)=k-1$. Then there exists a vertex of degree
$k-1$, say $u_1$. Let $N_G(u_1)=\{u_2,u_3,\cdots,u_k\}$. Choose
$S=\{u_1\}\cup N_G(u_1)$. Clearly, there is no pedant $S$-Steiner tree
in $G$. Hence $\tau_k(G)=0$, as desired. Suppose that $G$ is a graph satisfying Condition $(3)$. For the vertex set $S\subseteq V(G)$, there is no pedant $S$-Steiner tree
in $G$, and hence $\tau_k(G)=0$.
\qed
\end{pf}

\section{Nordhaus-Guddum type result}

In this section, we study the Nordhaus-Gaddum
type relations for pedant-tree connectivity.

\begin{thm}\label{th8}
Let $k,n$ be two integers with $3\leq k\leq n$, and let $G$ be a
connected graph of order $n$. Then

$(1)$ $0\leq \tau_k(G)+\tau_k(\bar{G})\leq n-k$;

$(2)$ $0\leq \tau_k(G)\cdot \tau_k(\bar{G})\leq [\frac{n-k}{2}]^2$.

Moreover, the upper and lower bounds are sharp.
\end{thm}
\begin{pf}
$(1)$ To avoid confusion, we denote the local pedant tree-connectivity of a $k$-subset
$S$ in a graph $G$ by $\tau(G;S)$. Since $G\cup \bar{G}=K_n$, for any $k$-subset $S$ we have
$\tau(G;S)+\tau(\bar{G};S)\leq \tau(K_n;S)$. Suppose that $\tau_k(K_n)=\tau(K_n;S_0)$ for some $k$-subset
$S_0$. Then we have
$$
\tau_k(K_n)=\tau(K_n;S_0)\geq \tau(G;S_0)+\tau(\bar{G};S_0)\geq \tau_k(G)+\tau_k(\bar{G}).
$$
This together with $\tau_k(K_n)= n-k$
results in $\tau_k(G)+\tau_k(\bar{G})\leq n-k$.

$(2)$ It follows immediately from $(1)$.
\end{pf}

\noindent {\bf Example 1:} Let $G'$ be a graph of order
$n-4$, and let $v_1v_2v_3v_4$ be a path. Let $G$ be the graph obtained
from $G'$ and the path by adding edges between the vertex $v_1$ and
all vertices of $G'$ and adding edges between the vertex $v_4$ and all
vertices of $G'$. Since $\delta(G)=\delta(\bar{G})=2$, it follows that $\tau_k(G)=\tau_k(\bar{G})=0$. So the lower bound of Theorem \ref{th8} is sharp for $3\leq k\leq n$. From Proposition \ref{pro2}, if $G=K_n$, then $\tau_k(G)=n-k$ and $\tau_k(\bar{G})=0$, and hence $\tau_k(G)+\tau_k(\bar{G})=n-k$. So the upper bound of Theorem \ref{th8} is sharp for $3\leq k\leq n$.

Let us focus on $(1)$ of Theorem \ref{th8}. If one of $G$ and
$\bar{G}$ is disconnected, we can characterize the graphs
attaining the upper bound by Lemma \ref{lem4}.

\begin{pro}\label{pro2}
For any graph $G$ of order $n$, if $G$ is disconnected, then
$\tau_k(G)+\tau_k(\bar{G})=n-k$ if
and only if $\bar{G}=K_n$.
\end{pro}

If both $G$ and $\bar{G}$ are all connected, we can obtain a
structural property of the graphs attaining the upper bound.

\begin{pro}\label{pro3}
If $\tau_k(G)+\tau_k(\bar{G})=n-k$,
then $\Delta(G)-\delta(G)\leq k-1$.
\end{pro}
\begin{pf}
Assume that $\Delta(G)-\delta(G)\geq k$. Since
$\tau_k(\bar{G})\leq \delta(\bar{G})=n-1-\Delta(G)$,
$\tau_k(G)+\tau_k(\bar{G})\leq
\delta(G)+n-1-\Delta(G)\leq n-1-k$, a
contradiction.\qed
\end{pf}\\

One can see that the graphs with
$\tau_k(G)+\tau_k(\bar{G})=n-k$
must have a uniform degree distribution.

From Corollary \ref{cor1}, the following observation are immediate.

\begin{obs}\label{obs3}
Let $G$ be a graph of order $n \ (n\geq 3)$. Then $\tau_n(G)+\tau_n(\bar{G})=0$ if and only if
$G$ is a graph of order $n$.
\end{obs}

From Corollary \ref{cor2} and Theorem \ref{th8}, we have the following result.

\begin{obs}\label{obs4}
Let $G$ be a graph of order $n \ (n\geq 4)$. Then
$$
0\leq \tau_{n-1}(G)+\tau_{n-1}(\bar{G})\leq 1.
$$
Furthermore, $\tau_{n-1}(G)+\tau_{n-1}(\bar{G})=1$ if and only if $G$ or $\bar{G}$ is complete; $\tau_{n-1}(G)+\tau_{n-1}(\bar{G})=0$ if and only if both $G$ and $\bar{G}$ are not complete.
\end{obs}

From Corollary \ref{cor3} and Theorem \ref{th8}, we have the following result.

\begin{pro}\label{pro4}
Let $G$ be a graph of order $n \ (n\geq 5)$. Then
$$
\tau_{n-2}(G)+\tau_{n-2}(\bar{G})=0~or~\tau_{n-2}(G)+\tau_{n-2}(\bar{G})=2.
$$
Furthermore, $\tau_{n-2}(G)+\tau_{n-2}(\bar{G})=2$ if and only if $G$ or $\bar{G}$ is complete; $\tau_{n-2}(G)+\tau_{n-2}(\bar{G})=0$ if and only if both $G$ and $\bar{G}$ are not complete.
\end{pro}
\begin{pf}
From Theorem \ref{th8},
$$
0\leq \tau_{n-2}(G)+\tau_{n-2}(\bar{G})\leq 2.
$$
Suppose $\tau_{n-2}(G)+\tau_{n-2}(\bar{G})=1$. Then $\tau_{n-2}(G)=1$ or $\tau_{n-2}(\bar{G})=1$. Without loss of generality, let $\tau_{n-2}(G)=1$ and $\tau_{n-1}(\bar{G})=0$. From Corollary \ref{cor3}, $G$ or $\bar{G}$ is a graph obtained from a complete graph of order $n$ by deleting at most two edges. Therefore, $e(G)={n\choose 2}-x$ and $e(\bar{G})={n\choose 2}-y$, where $1\leq x,y\leq 2$. Since $e(G)+e(\bar{G})={n\choose 2}$, it follows that $2\leq x+y={n\choose 2}\leq 4$ and hence $n=2$ or $n=3$, a contradiction. Hence
$$
\tau_{n-2}(G)+\tau_{n-2}(\bar{G})=0~or~\tau_{n-2}(G)+\tau_{n-2}(\bar{G})=2.
$$

Suppose $\tau_{n-2}(G)+\tau_{n-2}(\bar{G})=2$. Then $\tau_{n-2}(G)=2$ or $\tau_{n-2}(\bar{G})=2$. Without loss of generality, let $\tau_{n-2}(G)=2$ or $\tau_{n-2}(\bar{G})=0$. From Corollary \ref{cor3}, graph $G$ is complete. \qed
\end{pf}

\end{document}